# Pseudo Random test of prime numbers

Wang Liang[1]*, Huang Yan[2]

[1,2](Department of Control Science and Control Engineering, Huazhong University of Science and Technology, WuHan, 430074 P.R.China)

[2](Department of Mathematics, Huazhong University of Science and Technology, WuHan 430074 P.R.China)

**[Abstract]** The prime numbers look like a randomly chosen sequence of natural numbers, but there is still no strict theory to determine 'Randomness'. In these years, cryptography has developed a battery of statistical tests for randomness. In this paper, we just apply these methods to study the distribution of primes. Here the binary sequence constructed by second difference of primes is used as samples. We find this sequence can't reach all the 'random standard' of FIPS 140-1/2, but still show obvious random feature. The interesting self-similarity is also observed in this sequence. These results add the evidence that prime numbers is a chaos system.

**[Keywords]** Random test; Prime numbers

## 1 Introduction

Eratosthenes proved there are infinite prime numbers two thousand years ago. The prime number theorem describing general trend of prime's distribution is obtained one hundred year ago. But the detailed distribution of primes still keeps in secret. Recent research show prime numbers play an important role in physics and biology, so their detailed distribution has also come under increased scrutiny [1, 2].

When selecting a segment of prime numbers, they seem like the random collection in natural numbers. Applying statistical mechanics, M Wolf shows some interesting chaos pattern or said 'deterministic randomness' in prime numbers [3]. Some recent papers give the more detailed discussion for this random phenomenon [4-6]. Especially, the paper [7] suggests the distribution of difference of primes consistent with Gaussian statistics. A common character of these researches is the application of computer. The earlier works of probabilistic number theory are mainly theory research. A nice collection for them could be found in site [8].

"It is evident that the primes are randomly distributed but, unfortunately, we don't know what 'random' means." R.C. Vaughan (1990).

Obviously, defining whether a sequence is random is not an easy task. There are only few researches discussing the "random test" of primes. In a short paper [9], the author suggests the prime numbers are not randomly distributed by analyzing the distribution of first digital of primes. But this research only applies the Chi-square test. And first digital of primes has no much relation with the dynamical character of primes.

Recently, the need of high quality of random number producer in cryptography effectively prompts the randomness test research. Some standard random tests like FIPS have been developed [10]. We just apply these integrated tests to study the randomness of primes.

This paper is organized as follows. In the second section, how to construct the samples from primes is introduced. The following section show experimental results of FIPS random test. We also give the rough analysis for these results in this part. The last section is a short summary.

---
* Corresponding author: guoypm@hust.edu.cn

## 2 Samples data selection

Normally, the difference of consecutive prime numbers is used to study the detailed distribution of primes. We also apply this method.

Here $p_1, p_2, \cdots, p_i, p_{i+1}, \cdots$ is consecutive prime numbers, $p_i$ is ith prime number. Then their consecutive difference is: $C = \{c_i \mid c_i = p_{i+1} - p_i\}$. The second difference is

$$D = \{d_i \mid d_i = c_{i+1} - c_i\}, \quad i = 1, 2, 3, \cdots$$

We select the second difference of primes $D$ as basic test data. Its figure and statistic histogram is shown in fig.1. Here we calculate the consecutive prime numbers less than $1 \times 10^7$.

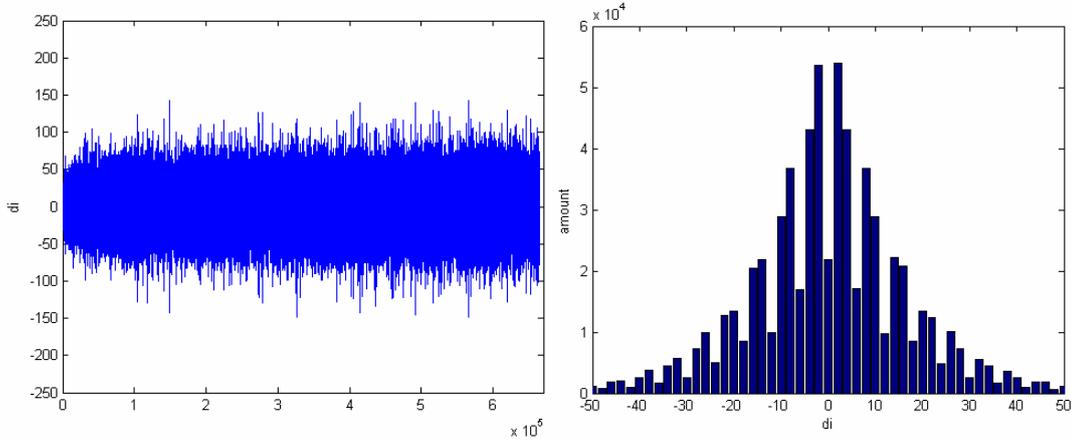

Fig 1.(a)The figure of $D$. (b) Statistic histogram statistic of $D$.

Then we convert series $D$ into a binary sequence $B$. For every element of $D$, the number smaller than 0 is remarked as 0. The number bigger than 0 is 1. Zero data is discarded. The mathematic expression for this process is:

$$B = (b_1, b_2, \cdots, b_i, \cdots) = (f(d_1), f(d_2), \cdots, f(d_i), \cdots), \quad b_i \in \{0, 1\}$$

Here, $f(d_i) = \begin{cases} 0, d_i < 0 \\ 1, d_i > 0 \\ discard, d_i = 0 \end{cases}$

The sequence $B$ is just our research object. There are three reasons that we select this binary sequence for random test. First, such symbolic description of data is a common method in statistics theory. Second, we could also find some supports from 'symbolic time series analysis' [11]. Related researches show such 'coarse description' could successfully feature the dynamical character of system [12]. And last, most random test like FIPS need the binary form of input data.

## 3 Random tests

Various statistical tests can be applied to a sequence to compare and evaluate the sequence to a truly random sequence. A statistical test is formulated to test a specific null hypothesis ($H_0$). Here null hypothesis under test is that the sequence being tested is random. Associated with this null hypothesis is the alternative hypothesis ($H_a$), which is that the sequence is not random. For each test, a relevant randomness statistic must be chosen and used to determine the acceptance ($H_0$) or rejection of the null hypothesis ($H_a$). The detailed description for random test could be found in [13].

We refer to the FIPS140-1/2 to design the random test for prime numbers. This test consists of four different tests that are applied to sequences of 20,000 bits.

Here we select two groups of data. The first group contains five consecutive 20000 bits sequence of $B$ beginning from the first prime number 2 (The biggest prime number less than $2 \times 10^6$). To make a clear comparison, another group comprising five same sequences starts from a prime number more than $1 \times 10^8$. The four FIPS random tests for these data are described respectively.

(1) Frequency Test (mono bit test)

The focus of the test is the proportion of zeroes and ones for the entire sequence. The purpose of this test is to determine whether the number of ones and zeros in a sequence are approximately the same as would be expected for a truly random sequence. FIPS request the number of ones must range in (9,725-10,275). To give a precise description for this result, we also apply the method of $\chi^2$ hypothesis test, here

$$\chi^2 = \frac{(n_0 - n_1)^2}{n}, n = n_0 - n_1,$$ $n_0$ is the number of '0' and the $n_1$ is the number of '1'.

If setting the significance level of 5%, we need $\chi^2 < 3.84$ to pass the random test.

The random test results for two groups of $B$ are shown in Table 1.

| G1 | $n_0$ | $n_1$ | $\chi^2$ | G2 | $n_0$ | $n_1$ | $\chi^2$ |
|---|---|---|---|---|---|---|---|
| 1 | 9976 | 10024 | 0.1152 | 1 | 9956 | 10044 | 0.3872 |
| 2 | 10029 | 9971 | 0.1682 | 2 | 10042 | 9958 | 0.3528 |
| 3 | 9981 | 10019 | 0.0722 | 3 | 10009 | 9991 | 0.0162 |
| 4 | 9973 | 10027 | 0.1458 | 4 | 9998 | 10002 | 0.0008 |
| 5 | 9976 | 10024 | 0.1152 | 5 | 9986 | 10014 | 0.0392 |

Table 1.(a) Statistical data of the first group. (b) The second group

We find two group sequence all can pass the frequency test. All subsequent tests depend on the passing of this test.

(2) Poker test

This test divides the 20000 bit stream into 5000 contiguous 4 bit segments. Count and store the number of occurrences of each of the 16 possible 4 bit values ('0000','0001',…,'1111'). Denote $f(i)$ as the number of each bit value where $1 \leq i \leq 15$. $f(i)$ has to satisfy following requirement:

$$\chi^2 = \frac{16}{5000}\sum_{i=0}^{15}[f(i)]^2 - 5000, \quad 2.16 < \chi^2 < 46.17$$

The $\chi^2$ of first group are (2173.3, 2.0069, 2153.0, 2178, 2201.1). The second group are (2069, 2217.1, 2171.2, 1951.9, 2245.4). Obviously, these values are much more than FIPS standard. This results shows proportions of different form of 4 bit block are not accordant. The figures $i \sim f(i)$ of two groups are shown in Fig 2. (X axis is i and Y corresponds to $f(i)$.)We could find an interesting result that the distributions of $f(i)$ of two groups are very similar. There is still no research for the random test of such sequence. So it's a new finding.

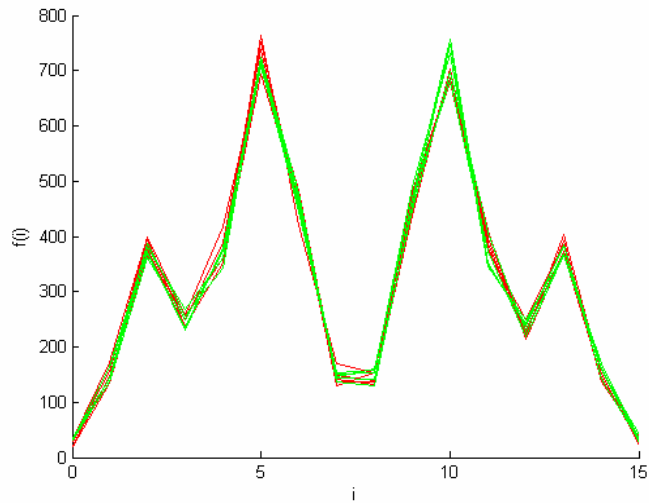

Fig 2. Distribution of $f(i)$. Red lines represent the data of the first group and the greens are the second group.

(3) Runs test

The focus of this test is the total number of runs in the sequence, where a run is defined as a maximal sequence of consecutive bits of either all ones or all zeros. The incidences of runs (for both consecutive zeros and consecutive ones) of all lengths ($\geq 1$) in the sample stream should be counted. For example, the '010101' contains the 6 one length of run (three '1' and three '0'), so its total runs are 6. The total runs of '000011111' are 2.

FIPS require the length of consecutive zeros and consecutive ones has to belong to well-defined ranges. Normally, the ideal total runs are the half length of sequence. Here the expected value is 10000. This test determines whether the oscillation between such zeros and ones is too fast or too slow.

Here the total runs of first group are (13660, 13570, 13691, 13664, 13694). The second groups are (13581, 13620, 13654, 13581, 13693). This result shows the oscillation in B is too fast. The distribution of different length run in one sequence is shown in Figure 3. (X axis is the length of run. The Y is its amount.) We could also find the distributions of run of two groups are very similar.

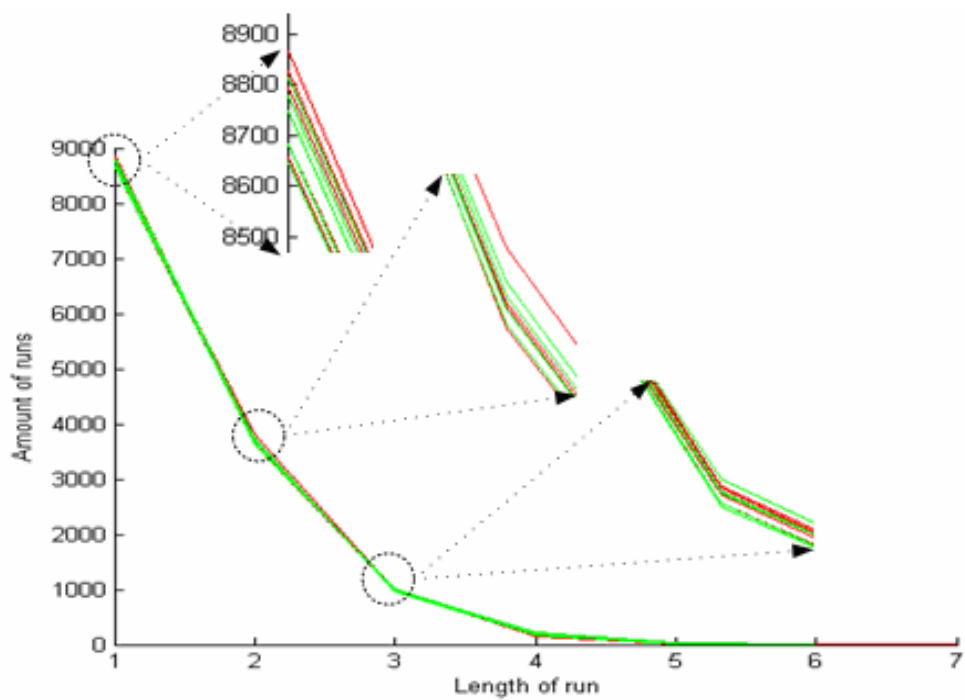

Fig3 Distribution of different length run. Red lines represent the data of the first group and the greens are the second group

(4) Long run test

FIPS needs no sequences of 26 consecutively equal bits occurring. For the first group, the longest runs are (6,6,6,6,7). The second group are (6,6,6,6,6). Not runes exceed 26.

We could draw the conclusion as follows: The binary sequence constructed from the second difference of primes can't pass all FIPS random test. Under this meaning, we could say the primes are not random distribution. But there are still obvious random feature in primes. In fact, the 'random' itself is not a strict concept. Moreover, the two groups of $B$ are selected from different segment of primes away from each other, but they have the similar statistical characters, which imply the self-similarity in primes. Here we still can't give the theory explanations for these findings in the detailed distribution of primes.

## 4 Summaries

Primes sequence shows obvious random feature but isn't the strict random sequence. It's the main results of this paper. We could safely say such 'deterministic randomness' is a kind of chaos phenomena.

The four tests of FIPS are only the basic random test. There have been nearly 200 different kinds of random tests. This paper shows the binary sequence obtained from the second difference of primes is appropriate to study the randomness of primes. Following this method, some 'universal law' of primes may be found in the future.